\documentclass{amsart}
\usepackage{amsfonts,amsmath,amssymb,amscd,latexsym,graphicx,epsfig,color,slashbox,a4wide}

\title[Fast Nielsen-Thurston classification of braids]{Fast Nielsen-Thurston classification of braids}

\theoremstyle{plain}
\newtheorem{thm}{Theorem}
\newtheorem{defi}[thm]{Definition}
\newtheorem{prop}[thm]{Proposition}
\newtheorem{lem}[thm]{Lemma}
\newtheorem{conj}[thm]{Conjecture}
\newtheorem{defprop}[thm]{Proposition-Definition}
\thanks{The author was partially supported by a grant from R\'{e}gion Bretagne, by MTM2010-19355 and FEDER, and by 
FONDECYT through postdoctoral grant no. 3130569.}

\begin{document}
\author{Matthieu Calvez}
\address{Matthieu Calvez, Departamento de Matem\'{a}tica y Ciencia de la Computaci\'{o}n, Facultad de Ciencias, Universidad de Santiago de Chile, Avenida Libertador Bernardo O'Higgins~3363, Santiago, Chile}
\email{calvez.matthieu@gmail.com}

\begin{abstract}
We prove the existence of an algorithm which solves
the reducibility problem in braid groups and runs in quadratic time with respect to the braid length for any fixed braid index.\end{abstract} 
%present a new algorithm for solving the reducibility problem in the braid groups. Its complexity is cubic with respect to braid length 
%for any fixed braid index. Our algorithm is based on a previous one by Gonz\'{a}lez-Meneses and Wiest. The latter was shown to have
%polynomial complexity (with respect to both braid length and braid index), modulo a conjecture about the speed of convergence of the cyclic sliding operation. 
%Thanks to a general result by Masur and Minsky about 
%mapping class groups, we partially establish the truth of this conjecture in the special case of rigid braids. Then, since any pseudo-Anosov braid has a small power conjugate to a rigid braid, we show that restricting our attention to this special rigid case is actually sufficient to get the above mentioned algorithm.\end{abstract}

%.
%\newpage

\maketitle
\section{Introduction} %and statements of the main results}
One of the main algorithmic decision problems regarding braids (viewed as mapping classes of a punctured closed disk) is the 
%so-called \emph{reducibility problem}, i.e. 
problem to determine the 
Nielsen-Thurston type of a given braid: reducible, periodic or pseudo-Anosov \cite{Fathi},\cite{Casson},\cite{Farb}. This problem is called the \emph{reducibility} problem because it amounts to 
determining whether a given non-periodic braid is reducible or not, i.e. whether it is reducible or pseudo-Anosov. Indeed, the case of periodic braids can be easily discarded: a braid~$x$ is periodic if and only if its $n$th power or its $(n-1)$st power is a power of the half-twist~$\Delta$ (see~\cite{bggm3}); and this is easy to decide algorithmically.

To solve the reducibility problem, two kinds of techniques have been used and several algorithms have been written; however none of them works in polynomial time 
with respect to the braid length for the general braid group~$B_n$. 

Firstly, the Bestvina-Handel algorithm \cite{bh} uses the theory of train-tracks and it is valid for any mapping class group. Although this algorithm works fast in practice, its 
theoretical complexity remains unknown.

Secondly, following the ideas introduced by Benardete, Gutierrez and Nitecki \cite{bgn,bgn1}, connections between the reducibility problem and the Garside structures of braid groups have been used for detecting reducibility \cite{bgn,bgn1,leelee,GMWiest,CalvezWiest,Calvezbkl}. 
Our work fits in this approach and builds mainly on the last algorithm by Gonz\'{a}lez-Meneses and Wiest \cite{GMWiest}.
% (see Section \ref{S:Garside} for an account on all of this). 

However Garside tools are not the only ones needed in the paper: we bring into play a very deep property of 
%-nonexplicit- existence result on 
Mapping Class Groups: the so-called \emph{linearly bounded conjugator property} 
%first due to Masur and Minsky 
\cite{MasurMinsky,jintao}, see Theorem \ref{masur}. At this point we already warn the reader that the algorithm given in the paper is not well-defined (although it will be actually described)
because it rests on the above linear bound, which is not explicitly known. Therefore our main result is an existence result only.
%and generalized later by Jing Tao 
%%\cite{jintao}. %first proved by  %Masur and Minsky \cite{MasurMinsky} for pseudo-Anosov mapping classes and extended by Jin Tao \cite{jintao} to all elements of any Mapping Class Group. 
%Unfortunately, the linear bounds from \cite{MasurMinsky,jintao} are not explicitly known; as a consequence, 
%%Because this linear bound is not explicitly known, 
%the algorithm constructed in the present paper is not 
%%
%%The reader should be warned that although we actually describe an algorithm as promised in the abstract, it is 
%well-defined (although it will be explicitly described) and our main result is an existence result only. 

The latter can be stated as follows. For a braid $x$, let us denote by $|x|$ the minimal possible length of a word representing $x$ whose letters are positive permutation braids and their inverses (in other words, the letters are braids in which any pair of strands crosses at most once and all crossings have the same orientation).

We will prove:
\begin{thm}\label{main}
Let $n$ be a positive integer. There exists an algorithm which decides the Nielsen-Thurston type of any given braid $x$ with $n$ strands and runs in time $O(|x|^2)$.
\end{thm}

The paper is organized as follows. In Section \ref{S:Garside} we recall useful tools from Garside theory and give precise statements relating the latter and the reducibility problem; an actual description of the algorithm whose existence is stated in Theorem \ref{main} is also given.
The detailed proofs are deferred to Section 3. 

% (see \cite{dehornoy} for the definition of a Garside group, see \cite{garside}, \cite{em}, \cite{GebhardtGM} for an introduction to the classical Garside structure of $B_n$, and \cite{bkldual} for the dual structure). 
%The best known so far is the latest algorithm due to Gonz\'{a}lez-Meneses and Wiest \cite{GMWiest}. The complexity is polynomial with respect to both 
%braid length
%and braid index, modulo a conjecture regarding the speed of convergence of the cyclic sliding operation $\mathfrak s$~(\cite{GebhardtGM}). In the present paper, we will give a partial proof of this conjecture, roughly : 
%
%\begin{thm} 
%There exists a constant $C(n)$ 
%\end{thm}

\section{The reducibility problem and Garside theory}\label{S:Garside}
\subsection{Reminders on Garside theory}
We first recall some basic notions of Garside theory in the specific case of braid groups, with emphasis on the \emph{classical} structure; references are \cite{garside,em,thurston}. 
The reader is referred to
%, although the definitions and the results can be extended \emph{mutatis mutandis} to this framework (
\cite{DehParis,dehornoy,DDGM} for a general account on Garside groups. 

The \emph{classical Garside structure} of the braid group consists in the following 2-fold data:~
%of the pair $(B_n^{+},\Delta)$, where~
$B_n^{+}$ is the monoid whose elements are braids which can be expressed as words on the Artin generators $\sigma_i$ with only positive exponents 
and $\Delta$ is the so-called \emph{Garside element} or half-twist.

The relation $\preccurlyeq$ on $B_n$ defined by $x\preccurlyeq y$ if and only if $x^{-1}y\in B_n^{+}$ defines a partial order called \emph{prefix order}, which turns out to be a lattice order. 
%The braid group $B_n$ is described by the well-known Artin's presentation 
%$$\left<\sigma_1,\ldots,\sigma_{n-1}\left | \begin{cases} \sigma_i\sigma_{j}\sigma_i=\sigma_j\sigma_i\sigma_j & \text{if}\ \  |i-j|=1\\
%                                                                  \sigma_i\sigma_j=\sigma_j\sigma_i & \text{if}\ \  |i-j|\geqslant 2\end{cases}\right .\right>.$$
%                                                                  
%The elements of $B_n$ which can be expressed as words on the letters $\sigma_i$ with only positive exponents form the monoid of \emph{positive braids} which is usually denoted by $B_n^{+}$. The relation $\preccurlyeq$ on $B_n$ defined by $x\preccurlyeq y$ if and only if $x^{-1}y\in B_n^{+}$ defines a partial order called \emph{prefix order}, which turns out to be a lattice order. 
%
%A key-r\^{o}le is played by the so-called \emph{Garside element} or half-twist: 
%$$\Delta=(\sigma_1\ldots\sigma_{n-1})\ldots(\sigma_1\sigma_2)\sigma_1.$$ 
We will denote by $\tau$ the inner automorphism associated to $\Delta$: this is an involution which for each $i$, maps $\sigma_i$ to $\sigma_{n-i}$; 
actually the center of $B_n$ ($n\geqslant 3$) is the cyclic group generated by $\Delta^2$.
%It is known that the square of $\Delta$ generates the center of $B_n$ (for $n\geqslant 3$); therefore the inner automorphism associated to $\Delta$ is an involution, which we denote by $\tau$ (note that it maps $\sigma_i$ to $\sigma_{n-i}$, for each~$i$).
 %We use the set of simple braids as a distinguished generating set of the braid group.
It can be shown that for any braid~$x$, there exist relative integers $r,s$ such that $\Delta^r\preccurlyeq x\preccurlyeq \Delta^s$. This allows to define the so-called \emph{infimum} and \emph{supremum} of $x$, respectively: 
%denoted by $\inf(x)$ and $\sup(x)$ 
$$\inf(x)=\max\{r\in \mathbb Z | \Delta^r\preccurlyeq x\}, \ \ \ \ \ \sup(x)=\min\{s\in\mathbb Z| x\preccurlyeq \Delta^s\}.$$ 
The \emph{canonical length} of $x$ is defined by $\ell(x)=\sup(x)-\inf(x)$. 

A central property of Garside groups is the existence of a distinguished generating set allowing for the definition of normal forms. 
Consider the set of positive prefixes of $\Delta$; these elements are called \emph{simple elements} or \emph{positive permutation braids} (because they are in one-to-one correspondence with the elements of the symmetric group on $n$ objects). Geometrically, simple elements are positive braids in which every pair of strands crosses at most once. Because it contains all Artin's generators $\sigma_i$, the set of simple elements generates $B_n$.

\begin{defi} 
A pair $(s_1,s_2)$ of simple elements is said to be \emph{left-weighted} if for any non-trivial positive prefix $t$ of $s_2$, the product $s_1t$ is not a simple element. 
\end{defi}

This allows to state:
\begin{prop} \cite{em}
Let $x\in B_n$. There exists a unique decomposition $x=\Delta^px_1\ldots x_r$, where $p=\inf(x)$, the $x_i$ are simple elements with $x_r\neq 1$ and (if $r\geqslant 2$) the pair $x_ix_{i+1}$ is left-weighted for $i=1,\ldots,r-1$. We have $\sup(x)=p+r$ and $\ell(x)=r$.
\end{prop}

Recall the braid length $|\cdot|$ defined in the introduction. It can be shown \cite{thurston} that every braid admits a unique decomposition of the form $a^{-1}b$, where $a,b$ are positive braids having no common non-trivial positive prefix. This is called the \emph{mixed canonical form}. Moreover, if $a=a_1\ldots a_k$ and $b=b_1\ldots b_l$ are the left normal forms of $a$ and $b$ respectively, it is shown in \cite{CharneyMeier} that the word $a_k^{-1}\ldots a_1^{-1}b_1\ldots b_l$ is a geodesic in the Cayley graph of $B_n$ with respect to the set of simple elements. Hence the braid length $|\cdot|$ is given by the length of the mixed canonical form. Finally notice that the latter is related to the canonical length in the following way:
if $x=\Delta^px_1\ldots x_r$ is the left normal form of $x\in B_n$, we have $|x|=\max(|p|,r,p+r)$ and $\ell(x)\leqslant |x|$.

Although we do not need that, it is worth mentionning that the braid group admits another Garside structure, called the \emph{dual Garside structure}, see 
 %and due to Birman, Ko and Lee 
 \cite{bkldual}. 
%The monoid involved is now the monoid whose elements are braids which can be written as positive words on the alphabet 
%$\{a_{i,j}=\sigma_,$
%the Garside element is $\delta=\sigma_1\ldots \sigma_{n-1}$.

The existence of normal forms for braids allows to construct algorithms for solving the conjugacy problem in the braid groups, that is for deciding whether any two given braids are conjugate and finding a conjugator whenever there exists one \cite{garside,em,Gebhardt,GebhardtGM,GebhardtGMalgorithm}. We now recall some related notions.

All existing algorithms for solving the conjugacy problem in braid groups rely on the definition of a particular finite computable characteristic subset of each conjugacy class in~$B_n$, consisting of its ``simplest" elements (in some sense).
A first example of such a characteristic subset is given by the following:
%
%that is an algorithm able to determine whether two given braids are conjugate or not and to find a conjugator whenever exists one. 
%This is usually done by constructing, in each conjugacy class, a finite subset which is characteristic of the class and consists of its "simplest" elements (in some sense).
\begin{defprop}\cite{em}
Let $x\in B_n$. The subset of the conjugacy class of~$x$ consisting of all elements with minimal canonical length is finite (and non-empty). 
Its elements have simultaneously maximal infimum and minimal supremum within the conjugacy class of $x$. 
This set is called the Super Summit Set of $x$ and is denoted by $SSS(x)$. 
\end{defprop}

An element in the Super Summit Set of a given braid can be computed using a special kind of conjugation called \emph{cyclic sliding}:
\begin{defi} \cite{GebhardtGM} 
Let $x=\Delta^px_1\ldots x_r$ be the normal form of $x\in B_n$. Suppose that $r>0$. The \emph{preferred prefix} of $x$ is the maximal positive prefix $t$ of $\tau^{-p}(x_1)$ such that $x_rt$ is a simple element; it is denoted by $\mathfrak p(x)$. If $r=0$, that is if $x$ is a power of $\Delta$, $\mathfrak p(x)$ is just defined to be the trivial braid.
The result of the conjugation of $x$ by its preferred prefix $\mathfrak p(x)$ is called \emph{cyclic sliding} of $x$ and denoted by $\mathfrak s(x)$. 
\end{defi}

The cyclic sliding operation actually achieves computing elements in Super Summit Sets, in an effective way. Indeed, only a polynomial (with respect to both length and braid index) number of iterations of $\mathfrak s$ is required to compute an element in the Super Summit Set: 
\begin{thm}\cite{bklbound,GebhardtGM}\label{T:BKLbound}
Let $x\in B_n$. The following alternative holds: $\ell(\mathfrak s^{\frac{n(n-1)}{2}-1}(x))<\ell(x)$ or $x\in SSS(x)$.
%%
%For any integer $t$ such that $t\geqslant \frac{n(n-1)}{2}\ell(x)$, we have $\mathfrak s^t(x)\in SSS(x)$.
\end{thm}

Another (smaller) example of conjugacy invariant subset was defined later by Gebhardt and Gonz\'{a}lez-Meneses:
\begin{defprop}\cite{GebhardtGM}
Let $x\in B_n$. The subset of the conjugacy class of $x$ consisting of all elements which are periodic points under the cyclic sliding operation is
finite and non-empty. It is a subset of the Super Summit Set of $x$, called set of Sliding Circuits of $x$ and denoted by $SC(x)$. 
\end{defprop} 
Iterative application of the cyclic sliding operation to a braid $x$ eventually reaches an element of $SC(x)$; however by contrast with the case of $SSS(x)$, no 
general bound on the number of repetitions involved is known (see Conjecture \ref{conj}).

To conclude this paragraph, we recall the important notion of \emph{rigid} braid: 
\begin{defi}\cite{GebhardtGM}
A braid $x$ is said to be \emph{rigid} if its preferred prefix is trivial. 
\end{defi}
In particular a rigid braid is a fixed point of $\mathfrak s$ and is an element of its own set of Sliding Circuits. Moreover, if a conjugacy class contains one rigid braid, then the corresponding set of Sliding Circuits consists only of rigid braids \cite{GebhardtGM}.

\subsection{Previous works}
We now turn to the relations between the reducibility problem and the notions above. First of all, we recall some definitions and properties related to reducible braids. 

Let $\mathbb D_n$ be the closed disk in $\mathbb C$ with diameter $[0,n+1]$ and with the points $\{1,\ldots,n\}$ removed. 
%(\{bgn}, \cite{leelee}, 
%\cite{CalvezWiest}, \cite{Calvezbkl}). 
As the Mapping Class Group of $\mathbb D_n$, the braid group $B_n$ induces a (right)-action on the set of isotopy classes of simple closed curves in $\mathbb D_n$. A curve is said to be \emph{non-degenerate} if it is not contractible and surrounds more than one and less than $n$ punctures. By the word ``curve" alone we will mean the isotopy class of a non-degenerate simple closed curve. 

A braid $x$ is \emph{reducible} if it preserves setwise a family of curves; such a curve is then said to be a \emph{reducing} curve for $x$. 
Notice that reducibility is a conjugacy-invariant property.
Moreover a reducible braid~$x$ is not periodic if the set of its \emph{essential} reducing curves, that is reducing curves which do not intersect any other reducing curve, is not empty. An important matter when we want to detect the reducibility of a braid is to actually detect reducing curves whenever they exist. This leads to formulate:

\begin{defi}\cite{bgn1,GMWiest}
A curve $\mathcal C$ in $\mathbb D_n$ is said to be \emph{round} if it is isotopic to a geometric circle in~$\mathbb D_n$. 
A curve $\mathcal C$ in $\mathbb D_n$ is said to be \emph{almost round} if it has a unique maximum and a unique minimum in the horizontal direction. 
The latter is equivalent to saying that $\mathcal C$ can be transformed into a round curve by the action of a simple element. 
\end{defi}

Whenever a reducible braid preserves a family of round curves, the reducibility is easy to detect (there are only $\frac{n(n-1)}{2}$ round curves), 
see \cite{bgn}. 
Moreover the notion of roundness of a curve behaves well with respect to Garside-theoretic operations:
\begin{thm}\cite{bgn1,juanjaca}\label{T:CyclicSlidingRound}
Let $x\in B_n$ be a reducible braid preserving a family of round curves. Then so does $\mathfrak s(x)$. 
\end{thm}

Notice that any family of curves can be transformed into a family of round curves under the action of some braid, in other word a reducible braid 
always preserves a family of round curves up to conjugacy. It follows from Theorems \ref{T:BKLbound} and \ref{T:CyclicSlidingRound} that for any reducible braid $x$, there exists some element of $SSS(x)$ that preserves a family of round curves. 

Therefore computing the whole Super Summit Set of a braid gives a manner of testing its reducibility: $x$ is reducible if and only if some element in $SSS(x)$ preserves a family of round curves.
We remark that the same approach can be carried out in the framework of the dual Garside structure, replacing round curves by \emph{standard} curves (see \cite{Calvezbkl}). 
However in either structure, the resulting algorithm is far from polynomial because the Super Summit Sets are exponentially large in general, with respect to both length and braid index \cite{juanjaca, prasolov}. 

In the special case of the four-strand braid group $B_4$, it is possible to overcome the latter difficulty, showing that \emph{every} element of the (classical) super summit set of a given reducible 4-braid with a reducing curve surrounding three punctures preserves a round or an almost-round curve. 
%it was shown in \cite{CalvezWiest} that there exists a polynomial-time algorithm 
This leads to a polynomial-time algorithm for solving the reducibility problem in $B_4$, described in \cite{CalvezWiest}. 
%
%Benardete, Gutierrez and Nitecki initiated the Garside-theoretical approach to the reducibility problem for braids in \cite{bgn} and \cite{bgn1}. They showed that the reducibility of a braid~$x$ is detectable 
%if one knows all the elements of the super summit set of $x$ (denoted by $SSS(x)$). 
%Indeed, it suffices to check whether some element of the super summit set preserves a family of \emph{round} curves. 
%The same can be done with the dual Garside structure and \emph{standard} curves replacing round curves~\cite{Calvezbkl}.
%Unfortunately, both algorithms are exponential because they demand computing the whole super summit sets, which in both structures have in general exponential size with respect to both braid length and braid index (see \cite{prasolov}, \cite{juanjaca}).

In \cite{leelee}, Lee and Lee replaced the condition about \emph{some} element of the \emph{super} summit set of a reducible braid $x$ by the condition that \emph{every} element of the \emph{ultra} summit set of~$x$ (which is another conjugacy invariant subset introduced in \cite{Gebhardt}) preserves a family of round curves. However this was shown at the cost of a technical hypothesis about the external and internal components of~$x$. 

In \cite{GMWiest}, Gonz\'{a}lez-Meneses and Wiest showed that every element of (some refined version of) the set of Sliding Circuits of a reducible braid preserves a family of round or almost-round curves. 

Both approaches from \cite{leelee} and \cite{GMWiest} suffer the same drawback, namely the lack of control on the distance to the first repetition when applying iteratively the cyclic sliding operation. The algorithm in \cite{GMWiest} indeed solves the reducibility problem in arbitrary braid groups with polynomial complexity, both in braid length
and braid index, provided the following conjecture holds:
%regarding the speed of convergence of the cyclic sliding operation $\mathfrak s$:

\begin{conj}\label{conj}{\rm{(}}\cite{GMWiest}, {\rm{Conjecture 3.5)}}
Let $x\in B_n$, with canonical length $r$. Let $t$ be the minimal positive integer such that $\mathfrak s^k(x)=\mathfrak s^t(x)$ for some $k$ with $0\leqslant k<t$. Then $t$ is bounded by a polynomial in $r$ and $n$.
\end{conj}

\subsection{Our algorithm}

One of the keys for proving Theorem \ref{main} will be a partial demonstration of Conjecture \ref{conj}. 
We shall prove, using Masur-Minsky's linear bound:
\begin{thm}\label{alternativepA} There exists a constant $C(n)$ (depending only on $n$) such that 
for any pseudo-Anosov $n$-strand braid $x\in SSS(x)$, the following holds:  
$x$ has a rigid conjugate if and only if $\mathfrak s^{C(n)\cdot |x|}(x)$ is rigid.
\end{thm}
%(see Theorem \ref{alternativepA} below). 

%We warn the reader that Theorem \ref{main} is only an existence result. We will actually describe an algorithm which is not well-defined because we do not know explicitly the constant $C(n)$ in 
%Theorem~\ref{alternativepA}, which comes from Masur-Minsky's conjugacy bound \cite{MasurMinsky}.
%Nevertheless, the existence of this constant is one of the keys for proving Theorem~\ref{main}:
%a partial proof of Conjecture 3.5 in~\cite{GMWiest} in the special case of rigid braids. Namely, we will prove: 

Theorem \ref{alternativepA} gives a partial solution to Conjecture \ref{conj}:
in the pseudo-Anosov rigid case, 
starting from a super summit element, Theorem \ref{alternativepA} guarantees that a rigid conjugate (or equivalently an element of the set of Sliding Circuits) is found  after only $C(n)\cdot|x|$ iterations of cyclic sliding (in other words, if a pseudo-Anosov super summit braid has rigid conjugates, then the cyclic sliding operation converges towards one of them in linear time with respect to braid length). This gives in this particular case (with the notation of Conjecture \ref{conj}) the bound $t\leqslant C(n)r+1$.

The importance of the rigid case comes from the following result, which will play a crucial role in our proof of Theorem \ref{main}. It is due to Birman, Gebhardt and Gonz\'{a}lez-Meneses:

\begin{thm}\cite{bggm1}.\label{bggm1}
Let $x\in B_n$ be a pseudo-Anosov braid. There exists a positive integer $m<{(\frac{n(n-1)}{2}})^3$ such that $x^m$ is conjugate to a rigid braid. 
\end{thm}
%Namely, we will see that if $x\in SSS(x)$ has a rigid conjugate, then such a rigid conjugate can be obtained from $x$ by applying $k$ iterations of cyclic sliding and we will show that $k$ is bounded above linearly by $|x|$ for any fixed number of strands. 
 Finally we recall the two following results from \cite{GMWiest}:

\begin{thm}{\rm(\cite{GMWiest}, Theorem 5.16)}.\label{5.16}
Let $x\in B_n$ be a non-periodic, reducible braid which is rigid.
Then some essential reducing curve of $x$ is round or almost-round. More precisely, 
there is some positive integer $k\leqslant n$ such that one of the following holds:
\begin{itemize}
\item[(1)] $x^k$ preserves a round essential curve, 
\item[(2)] $\inf(x^k)$ and $\sup(x^k)$ are even and either $\Delta^{-\inf(x^k)}x^k$ or $x^{-k}\Delta^{\sup(x^k)}$ is a positive braid preserving an almost-round essential reducing curve whose interior strands do not cross.
\end{itemize}
\end{thm}
\begin{thm}{\rm(\cite{GMWiest}, Theorem 2.9)}.\label{2.9}
There is an algorithm which decides whether a given positive braid $x$ preserves an almost-round curve whose interior strands do not cross. Its complexity is $O(\ell(x)n^4)$.
\end{thm}

We are now ready to describe the algorithm promised in Theorem \ref{main}. It takes as input an $n$-braid~$x$. The output is ``periodic", ``reducible" or ``pseudo-Anosov".
\begin{itemize}
\item[1.] If $x^{n-1}$ or $x^n$ is a power of $\Delta$, return ``periodic" and stop.
\item[2.] For $i=1,\ldots,(\frac{n(n-1)}{2})^3-1$ compute the normal form of $x^i$. Iteratively apply cyclic sliding to~$x^i$ until the canonical length has not decreased during the $\frac{n(n-1)}{2}-1$ last iterations. This computes $y_i\in SSS(x^i)$. Then compute $z_i=\mathfrak s^{C(n)\cdot|y_i|}(y_i)$. If none of the braids $z_i$ is rigid return ``reducible" and stop. Else let $j$ be such that~$z_j$ is rigid.
%\item[3.] If $z_j$ preserves a family of round curves, return ``reducible" and stop.
\item[3.] For $k=1,\ldots,n$, apply the algorithm in \cite{bgn} to the braid $z_j^k$ to test whether it preserves a round curve; apply the algorithm in Theorem \ref{2.9} to both braids $\Delta^{-\inf(z_j^k)}z_j^k$ and $z_j^{-k}\Delta^{\sup(z_j^k)}$. If a round or an almost-round reducing curve is found, then return ``reducible" and stop.
\item[4.] Return ``pseudo-Anosov".
\end{itemize}
As mentioned above, we remark that this algorithm, and specifically Step 2, is not well-defined because the constant $C(n)$ is not explicitly known. In the next section, we will prove Theorem \ref{alternativepA}, show the correctness of the above algorithm and study its complexity.

\section{Proofs of our results}
In order to prove Theorem \ref{alternativepA}, we combine two important results. The first one is the already mentionned linearly bounded conjugator property for Mapping Class Groups \cite{MasurMinsky}. 
Although the range of surfaces considered by Masur and Minsky is much broader, only the ($n+1)$-times punctured sphere $\mathbb S_{n+1}$ ($n\geqslant 3$) is relevant for our purposes, so we state their result in this special case: 
%of the $(n+1)$-times punctured sphere $\mathbb S_{n+1}$.

\begin{thm}{\rm{(\cite{MasurMinsky}, Theorem 7.2)}}.\label{masur} %Let $n\in \mathbb N$. Let $\mathcal S$ be any system of generators for $\mathcal {MCG}(\mathbb D_n)$ and let $|\cdot|$ be the length of 
%mapping classes with respect to $\mathcal S$. 
Let~$\mathcal G$ be any generating set of the Maping Class Group $\mathcal {MCG}(\mathbb S_{n+1})$. 
There exists a constant~$\gamma(\mathcal G)$, depending only on $\mathcal G$, such that any pair of conjugate pseudo-Anosov mapping classes
%$n$-strand braids $x,y$ 
can be related by a conjugating element~$w$ satisfying
$$|w|_{\mathcal G}\leqslant \gamma(\mathcal G)(|x|_{\mathcal G}+|y|_{\mathcal G}),$$
(where $|\cdot|_{\mathcal G}$ denotes the word length with respect to the chosen generating set~$\mathcal G$). 
\end{thm}

We aim at an analogous result for braids, namely we want to show:

\begin{prop}\label{masurbraids} %Let $n\in \mathbb N$. Let $\mathcal S$ be any system of generators for $\mathcal {MCG}(\mathbb D_n)$ and let $|\cdot|$ be the length of 
%mapping classes with respect to $\mathcal S$. 
There exists a constant~$c(n)$, depending only on $n$, such that any pair of conjugate pseudo-Anosov $n$-braids
%$n$-strand braids $x,y$ 
can be related by a conjugating element~$w$ satisfying
$$|w|\leqslant c(n)(|x|+|y|).$$
\end{prop}

%We notice that the exact statement of Theorem 7.2 in \cite{MasurMinsky} concerns mapping class groups of closed surfaces with punctures, in particular the mapping class group of the $(n+1)$-times punctured sphere $\mathcal{MCG}(\mathbb S_{n+1})$, not the braid group. 
%\begin{proof} 

Before proving Proposition \ref{masurbraids}, we recall that the quotient $B_n/\left<\Delta^2\right>$ can be seen as the Mapping Class Group of an $n$-times punctured closed disk (with boundary fixed setwise). For a braid $x$ in $B_n$, denote by
$\hat x$ its image in the quotient $B_n/\left<\Delta^2\right>$. Simple elements are sent bijectively to a generating set of the quotient (whose elements we still call simple); this defines a length $||\cdot||$ on $B_n/\left<\Delta^2\right>$
%The set of simple elements induces a word length on the quotient $B_n/\left<\Delta^2\right>$ which we denote by~$||\cdot||$ 
(notice that for any braid~$x$, 
$||\hat x||\leqslant |x|$).
%${||x \pmod {\Delta^2}||\leqslant |x|}$).
Collapsing the boundary of the $n$-times punctured closed disk to a puncture in the sphere~$\mathbb S_{n+1}$ (where the punctures are uniformly placed along the horizontal great circle), we can view $B_n/\left<\Delta^2\right>$ as the finite index subgroup of $\mathcal{MCG}(\mathbb S_{n+1})$ consisting 
of the mapping classes which fix the $(n+1)$st puncture.
%identify the quotient $B_n/\left<\Delta^2\right>$ with the finite index subgroup of $\mathcal {MCG}(\mathbb S_{n+1})$ consisting of the mapping classes which fix the $(n+1)$st puncture (i. e. the puncture corresponding to the boundary of the~$n$-times punctured closed disk after collapsing the boundary into a puncture in the sphere $\mathbb S_{n+1}$).
The group $\mathcal{MCG}(\mathbb S_{n+1})$ is equipped with the generating set
%$$\mathcal G_n=\{\hat \sigma_1,\ldots,\hat \sigma_{n-1},\rho\},$$
%where $\rho$ is a clockwise rotation by an angle of~$\frac{2\pi}{n+1}$. 
consisting of the simple elements in the quotient $B_n/\left<\Delta^2\right>$ together with a clockwise rotation by an angle of~$\frac{2\pi}{n+1}$, which we denote by $\rho$. 
Notice that for any $u\in B_n/\left<\Delta^2\right>$, we have $|u|_{\mathcal G_n}\leqslant ||u||$. Conversely, we will see in a simple computational way that the length $||u||$ is linearly bounded in terms of $|u|_{\mathcal G_n}$: 
\begin{lem}\label{calcul}
For any $u\in B_n/\left<\Delta^2\right>$, we have $||u||\leqslant \frac{n(n-1)}{2} |u|_{\mathcal G_n}$. 
\end{lem}
\begin{proof}
Let $u\in B_n/\left<\Delta^2\right>$. We will construct a word representative $W$ of $u$ using only the letters $\hat\sigma_1^{\pm1},\ldots,\hat\sigma_{n-1}^{\pm1}$ such that $||W||\leqslant \frac{n(n-1)}{2} |u|_{\mathcal G_n}$. 
This is achieved thanks to the following relations in $\mathcal{MCG}(\mathbb S_{n+1})$ (which can be easily deduced from the presentation 
in \cite{birman} Theorem 4.5).
For $1\leqslant i\leqslant n-1$, $0\leqslant j\leqslant n$: 
\begin{eqnarray}
\hat\sigma_i\rho^j & = & \rho^j\hat\sigma_{i-j\pmod {n+1}},\ \ \ i-j\neq 0,n\pmod {n+1}\\
\hat\sigma_i\rho^i & = & \rho^{i+1}(\hat\sigma_1...\hat\sigma_{n-1})^{-1},\\
\hat\sigma_i\rho^{i+1}& = & \rho^{i}(\hat\sigma_{n-1}\ldots \hat\sigma_1)^{-1}.
\end{eqnarray}
%Set $$\mathcal U=\{\rho^{-j}\hat\sigma_i\rho^j,\  i=1,\ldots,n-1,\ j=0,\ldots,n,\ i-j\neq 0,n\pmod{n+1} \cup \{\rho^{-i-1}\hat\sigma_i\rho^i,\rho^{-i}\hat\sigma_i\rho^{i+1},\  i=1,\ldots,n-1\}.$$ 
%It is easily seen that $\mathcal U$ is a generating set of $B_n/\left<\Delta^2\right>$, and that $|u|_{\mathcal U}\leqslant |u|_{\mathcal G_n}$, for any $u\in B_n/\left<\Delta^2\right>$. 
This allows to gather (at the beginning) all powers of $\rho$ appearing in a shortest representative for~$u$ on the alphabet $\mathcal G_n$; this results in a power of $\rho^{n+1}$ (because $u\in B_n/\left<\Delta^2\right>$) followed by a word $W$ in $\hat\sigma_1^{\pm1},\ldots,\hat\sigma_{n-1}^{\pm1}$. The word $W$  represents $u$ and 
as a word on the simple elements and their inverses, its length is bounded by $\frac{n(n-1)}{2}|u|_{\mathcal G_n}$ (because each simple element can be written with at most $\frac{n(n-1)}{2}$ letters~$\hat\sigma_i$ so that the above relations are used at most $\frac{n(n-1)}{2}|u|_{\mathcal G_n}$ times; and because both $\hat\sigma_1\ldots\hat\sigma_{n-1}$ and $\hat\sigma_{n-1}\ldots\hat\sigma_1$ are simple). 
%Indeed, the number of letters 
%$\sigma_i$ ($i=1,\ldots,n-1$) concerned by the transformation is certainly bounded by $\frac{n(n-1)}{2}|u|_{\mathcal G_n}$ (because each simple element can be written with at most $\frac{n(n-1)}{2}$ letters~$\sigma_i$) and we see that each of them gives rise to a single simple element in the new obtained word. 
\end{proof}

{\it{Proof of Proposition \ref{masurbraids}.}}
Given a pair of conjugate pseudo-Anosov $n$-braids $x$ and~$y$, we know a conjugating element between~$\hat x$ and $\hat y$, say $\upsilon$ in the quotient $B_n/\left<\Delta^2\right>$. 
In their proof of Theorem~\ref{masur}, Masur and Minsky construct a ``short" conjugating element $\upsilon'$ between~$\hat x$ and $\hat y$; this element~$\upsilon'$ is expressed as a product 
%Masur-Minsky's proof of Theorem \ref{masur} constructs a ``short" conjugating element $\upsilon'$ between~$\hat x$ and~$\hat y$ in $\mathcal {MCG}(\mathbb S_{n+1})$. Moreover, this element $\upsilon'$ is given as a product 
$\hat x^m\upsilon$, for some integer $m$ and 
%whose length is linearly bounded by the length of $x$ and $y$. 
%The mapping class $\upsilon'$ is the product $\hat x^m\upsilon$, for some integer~$m$ and 
therefore it belongs to the subgroup $B_n/\left<\Delta^2\right>$ of 
$\mathcal{MCG}(\mathbb S_{n+1})$. Moreover, $$|\upsilon'|_{\mathcal G_n}\leqslant \gamma(\mathcal G_n)(|\hat x |_{\mathcal G_n}+|\hat y|_{\mathcal G_n})\leqslant \gamma(\mathcal G_n)(||\hat x||+||\hat y||)$$ and we get from Lemma \ref{calcul}
$${||\upsilon'||\leqslant \frac{n(n-1)}{2}\gamma(\mathcal G_n)(||\hat x||+||\hat y||)}.$$

%mapping classes are supposed to fix  %ie is shown in the usual mapping class group setting, that is, mapping classes are assumed to fix 
%setwise the boundary of the disk (whereas braids fix it pointwise). 
%In other words, the mentioned result from \cite{MasurMinsky} actually holds in the quotient $B_n/\left<\Delta^2\right>$ (with the same set of generators). 

Finally, as $\left<\Delta^2\right>$ is the center of $B_n$, and because a braid $x$ conjugate to $y$ cannot be conjugate to
$\Delta^{2k}y$ for $k\neq 0$, any lifting of $\upsilon'$ in $B_n$ conjugates $x$ to $y$ and we can choose one, say $w$, so that $|w|=||\upsilon'||$. Therefore, taking $c(n)=\frac{n(n-1)}{2}\gamma(\mathcal G_n)$ 
achieves the proof of 
Proposition \ref{masurbraids}. \hfill$\Box$
%so that Masur-Minsky's result can actually be translated to braids, as stated in Proposition \ref{masurbraids}.

The second step towards Theorem \ref{alternativepA} is a general fact about Garside groups. It explains that if a super summit element 
has a rigid conjugate, then iterated cyclic sliding is the shortest way of obtaining such a rigid conjugate.
\begin{thm}\cite{GebhardtGM}.\label{gebhardtgm}
Let $x\in B_n$ and assume that $x$ is conjugate to a rigid braid.
\begin{itemize}
\item[(1)] There exists a unique positive braid $f(x)$ such that $f(x)^{-1}x{f(x)}$ is rigid and ${f(x)\preccurlyeq g}$ for any positive braid $g$ such that $g^{-1}xg$ is rigid. 
\item[(2)] If $y\in SSS(x)$, then (unless $y$ is already rigid) there exists some positive integer $k$ such that $f(y)=\prod_{i=1}^{k} \mathfrak p(\mathfrak s^{i-1}(y))$. That is, $f(y)$ is the product 
of the $k$ conjugating simple elements involved when applying $k$ iterations of cyclic sliding to $y$.
%$f_k(y)=\prod_{i=1}^{k} \mathfrak p(\mathfrak s^{i-1}(y))$ satisfies $y^{f_k(y)}$ is rigid. Moreover, $f(y)=f_k(y)$ is the minimal positive braid doing so,
%i.e. for any positive braid $g$ such that $y^g$ is rigid, we have $f(y)\leqslant g$.
\end{itemize}
\end{thm}

Now, the proof of Theorem \ref{alternativepA} is just a combination of both of the previous results.

\textit{Proof of Theorem \ref{alternativepA}.} Let $x$ be a pseudo-Anosov $n$-strand braid such that ${x\in SSS(x)}$. Let us assume that $x$ has a rigid conjugate $z$.
By Proposition \ref{masurbraids}, there exists $w\in B_n$ such that $z=w^{-1}xw$ and $|w|\leqslant c(n)(|x|+|z|)$. 
Since $x,z\in SSS(x)$, we have $|x|=|z|$.
%Let $m$ be the minimal positive integer such that $\delta^m$ is central ($m=2$ for the classical Garside structure, $m=n$ for the dual one).
Let $r$ be the number of negative factors in the mixed canonical form of~$w$. If $r$ is even, then $w'=\Delta^rw$ is a positive braid conjugating $x$ to $z$ 
(recall that $\Delta^2$ is central).
%is a multiple of $m$, then $w'=\delta^kw$ is a positive braid conjugating $x$ to $z$.
Otherwise $r$ is odd and $w'=\Delta^{r+1}w$ does the same job.
%we denote by $\bar k$ the rest of the division of $k$ by $m$ and the braid $w'=\delta^{k+m-\bar k}w$ is a positive braid conjugating $x$ to $z$. 
In either case, we get a positive braid $w'$ conjugating $x$ to $z$ with $|w'|\leqslant |w|+1\leqslant (2c(n)+1)|x|$ (we may assume that $|x|\geqslant 1$). 

Let $k$ and $f(x)=\prod_{i=1}^{k} \mathfrak p(\mathfrak s^{i-1}(x))$ be as in Theorem \ref{gebhardtgm}. Then $f(x)\preccurlyeq w'$. It follows that $|f(x)|\leqslant |w'|$. 
%Let $q=\frac{n(n-1)}{2}$ be the length of $\Delta$ with respect to the atoms. 
As the braid $f(x)$ is a product of $k$ simple elements, we have $k\leqslant \frac{n(n-1)}{2}|f(x)|$ (because a simple element can be written with at most $\frac{n(n-1)}{2}$ letters $\sigma_i$) so that finally $k\leqslant \frac{n(n-1)}{2}\cdot(2c(n)+1)|x|$. Thus taking $C(n)=\frac{n(n-1)}{2}\cdot(2c(n)+1)$ (which depends only on~$n$), we have shown the following:~$x$ has a rigid conjugate if and only if $\mathfrak s^{C(n)\cdot|x|}(x)$ is a rigid braid (notice that $\mathfrak s^m(z)=z$ for any rigid braid $z$ and any integer $m\in \mathbb N$).~\hfill$\Box$

We notice that Theorem \ref{alternativepA} can be shown as well in the dual setting but we will not need this. We now turn to the proof of Theorem \ref{main}, showing the correctness of the algorithm in Section \ref{S:Garside} and studying the complexity of each step.

{\it{Proof of Theorem \ref{main}.}} 
The correctness of Step 1 is shown in~\cite{bggm3}. This step just consists in a computation of left normal form; therefore it takes time $O(\ell(x)^2)$ for any fixed $n$, according to~\cite{thurston}.

Let us prove that Step 2 is correct. First, by Theorem \ref{T:BKLbound}, 
%%that starting from any $n$-braid $x$, either $x$ is an element of $SSS(x)$ or iteratively applying the cyclic 
%sliding operation $(\frac{n(n-1)}{2}-1)$ times to $x$ decreases the canonical length. Therefore for each $i$, 
the braid~$y_i$ is an element of $SSS(x^i)$ for each $i$.
Then if $x$ is a pseudo-Anosov braid,
by Theorem~\ref{bggm1}, at least one of the braids~$x^i$ (and therefore~$y_i$) is pseudo-Anosov with a rigid conjugate and by Theorem~\ref{alternativepA}
at least one of the braids~$z_i$ is rigid. 

Let us calculate the complexity of Step 2. The computations of normal forms are known to be quadratic with respect to the length \cite{thurston}. 
We then recall that each instance of cyclic sliding (when applied to a braid already in normal form) has linear complexity with respect to the braid length for any given braid index (see \cite{GebhardtGMalgorithm}, Theorem 4.4). 
Therefore for any $i=1,\ldots ,(\frac{n(n-1)}{2})^3-1$, the complexity of computing $y_i$ (which requires at most $(\ell(x^i)-1)\cdot(\frac{n(n-1)}{2}-1)$ iterations of cyclic sliding) is quadratic with respect to 
the braid length whenever $n$ is fixed. The same is true for the computation of~$z_i$ from $y_i$ which requires $C(n)|y_i|\leqslant C(n)i|x|$ iterations of cyclic sliding.
%$$O\left(\left((\frac{n(n-1)}{2}-1)(|x^i|-1)\right)|x^i|^2n^2+(C|y_i|)|y_i|^2n^2\right).$$
%This computation has complexity 
%$O(\ell(x)^2n^2)$. Therefore for any $i=1,\ldots ,(\frac{n(n-1)}{2})^3-1$, computing $y_i$ is doable in time $O(i^2|x|^2n^2)$. Then, 

%The correctness of Step~2 can be shown by contradiction: if $x$ was non-periodic pseudo-Anosov, 
%then by Theorem~\ref{bggm1}, some of the braids $x^i$ (and therefore~$y_i$) would be pseudo-Anosov with a rigid conjugate and by Theorem~\ref{alternativepA}
%at least one of the braids $z_i$ would be rigid. 
%To evaluate the complexity of that step, we need to recall:
%\begin{thm}\cite{bklbound}.There is an algorithm which takes as its input a braid $x$ with~$n$ strands, runs for time $O(\ell(x)^2n^2)$, and outputs some element $x'\in SSS(x)$.
%\end{thm}
% and the computation of $z_i$ requires $C |y_i|\leqslant Ci|x|$ iterations of $\mathfrak s$, 
%each of them being performed in time $O(|x|^2)$ for a given $n$ (. Finally Step 2 is doable in cubic time with respect to $|x|$.
 
The validity of Step 3 follows from Theorem \ref{5.16}. This step consists in applying the algorithm in~\cite{bgn} to at most $n$ braids of length at most $nj|x|$ and the algorithm of Theorem \ref{2.9}
to at most $2n$ braids of length at most $nj|x|$. Both of these algorithms work in linear time with respect to the length so that Step 3 is linear with respect to $|x|$.
%the correctness of the algorithm in~\cite{GMWiest}. 
%Moreover, Step 3 has complexity bounded above linearly by the length of $z_j$, hence by the length of $x$, following the well-known algorithm in \cite{bgn}. Finally, the validity of Step 4 follows from Theorem \ref{5.16}. This step consists in applying the algorithm of Theorem \ref{2.9} to $2n$ braids of length at most $nj|x|$.
%According to Theorem \ref{2.9} the complexity of that step is linear with respect to $|x|$.
%This achieves the proof of Theorem \ref{main}. 
\hfill$\Box$

We notice that the present algorithm does not always yield the knowledge of reducing curves for reducible elements (actually this failure happens when reducibility is detected at Step~2). Thus, in view of the program in \cite{bggm1}, \cite{bggm2},~\cite{bggm3}, writing an algorithm for solving the conjugacy problem in braid groups in polynomial time still requires the following:
\begin{itemize}
\item[(i)] find explicitly the constant $C(n)$ to make the algorithm in~Theorem \ref{main} explicit. This amounts to bounding explicitly the required number of cyclic slidings to obtain (if it exists) a rigid conjugate from a pseudo-Anosov super summit element (see Theorem \ref{alternativepA}); alternatively this rests on the knowledge of an explicit value for Masur-Minsky's constant~$c(n)$ (see Proposition \ref{masurbraids}),
\item[(ii)] find reducing curves of a braid, in polynomial time, whenever the braid is known to be reducible,
\item[(iii)] find a polynomial bound on the number of rigid braids in a given pseudo-Anosov conjugacy class.
\end{itemize}
%Moreover, if we want the algorithm to be explicitly given, we still need to determine explicitly the constant $C$ of Theorem \ref{alternativepA}. 

We finish with a discussion of the special case of the four-strand braid group $B_4$. 
If we want to decide the Nielsen-Thurston type of a given 4-braid, the algorithm in~\cite{CalvezWiest} should rather be used instead of the present one because it is implementable and 
it finds explicitly the reducing curves whenever they exist (in polynomial time).
%In this case, according to \cite{CalvezWiest}, we know how to do (i). 
Using the Birman-Ko-Lee structure of $B_4$, the author together with Bert Wiest show in \cite{CalvezWiest4Strand} the existence of a bound as in (iii) (which depends on 
Masur-Minsky's constant $c(4)$, see Proposition \ref{masurbraids}).
%In a forthcoming paper \cite{CalvezWiest4Strand}, the author together with Bert Wiest show in the case of $B_4$ the existence of a bound as in (ii). 
Unfortunately, they do not know yet how to make explicit the constant~$c(4)$ (nor $C(4)$), so that the cardinality of the ultra summit set of a pseudo-Anosov rigid 4-braid is not explicitly known. Nevertheless \cite{CalvezWiest4Strand} presents a polynomial-time algorithm for solving the conjugacy problem in $B_4$.
%According to \cite{CalvezWiest}, the first difficulty can be overcome in the special case of the 4-strand braid group~$B_4$. In a forthcoming paper, the author together with Bert Wiest show how to bypass the second of these difficulties in $B_4$; they can also explicit the constant $C$and therefore obtain a polynomial-time algorithm for solving the conjugacy problem in~$B_4$.

{\bf{Acknowledgements.}} 
The author thanks Jin Tao for her help with Proposition \ref{masurbraids} and Juan Gonz\'{a}lez-Meneses for useful comments on an earlier version of that work.
His gratefulness also goes to Bert Wiest for many inspiring discussions and careful readings of the manuscript. 
%He also would to thank Juan Gonz\'{a}lez-Meneses and Bert Wiest for useful comments on an earlier version of that work.


\begin{thebibliography}{1}
\bibitem{bgn} D. Benardete, M. Gutierrez, Z. Nitecki, \textit{A combinatorial approach to reducibility of mapping classes}, Mapping Class Groups and moduli spaces of Riemann surfaces 
(G\"{o}ttingen, 1991/Seattle, WA, 1991), 1-31, Contemp. Math., 150, Amer. Math. Soc., Providence, RI, 1993.
\bibitem{bgn1} D. Benardete, M. Gutierrez, Z. Nitecki, \textit{Braids and the Nielsen-Thurston classification}, J. Knot Theory Ramifications, 4 (1995), 549-618.
\bibitem{bh}M. Bestvina, M. Handel, \textit{Train-tracks for surface homeomorphisms}, Topology 34 (1995), no.~1, 109-140.
\bibitem{birman} J. Birman, \textit{Braids, Links and Mapping Class Groups}, Annals of Math. Studies 82, (1974).
\bibitem{bggm1} J. Birman,V. Gebhardt, J. Gonz\'{a}lez-Meneses, \textit{Conjugacy in Garside Groups I: Cycling, Powers and Rigidity}, Groups Geom. Dyn. 1 (2007), no. 3, 221-279.
\bibitem{bggm2} J. Birman,V. Gebhardt, J. Gonz\'{a}lez-Meneses, \textit{Conjugacy in Garside Groups II: Structure of the ultra summit set}, Groups Geom. Dyn. 2 (1), (2008), 16-31.
\bibitem{bggm3} J. Birman, V. Gebhardt, J. Gonz\'{a}lez-Meneses, \textit{Conjugacy in Garside groups III: Periodic braids}, J. Algebra 316 (2), (2007), 746-776.
\bibitem{bkldual} J. Birman, K.-H. Ko, S.-J. Lee, \textit{A new approach to the word and conjugacy problems in the braid groups}, Adv. Math. 139 (2) (1998), 322-353.
\bibitem{bklbound} J. Birman, K.-H. Ko, S.-J. Lee, \textit{The infimum, supremum, and geodesic length of a braid conjugacy class}, Adv. Math. 164 (2001), 41-56.
\bibitem{Calvezbkl} M. Calvez, \textit{Dual Garside structure and reducibility of braids}, J. Algebra 356 (1), (2012), 355-373.
\bibitem{CalvezWiest4Strand} M. Calvez, B. Wiest, \textit{The conjugacy problem in the four-strand braid group}, arXiv:1204.6507.
\bibitem{CalvezWiest} M. Calvez, B. Wiest, \textit{Fast algorithmic Nielsen-Thurston classification of four-strand braids}, J. Knot Theory Ramifications, 21 (5), (2012).  
\bibitem{Casson} A. Casson, S. Bleiler, \textit{Automorphisms of surfaces after Nielsen and Thurston}, LMS Student Texts, 9. Cambridge University Press, Cambridge, 1988. 
\bibitem{CharneyMeier} R. Charney, \textit{Geodesic automation and growth functions for Artin groups of finite type}, Math. Ann. 301 no.2, (1995), 307-324.
\bibitem{dehornoy} P. Dehornoy, \textit{Groupes de Garside}, Ann. Sci. \'{E}cole Norm. Sup. (4) 35 (2002), no. 2, 267-306. 
\bibitem{DDGM} P. Dehornoy, F. Digne, E. Godelle, D. Krammer, J. Michel, \textit{Foundations of Garside Theory}, in progress.
\bibitem{DehParis} P.Dehornoy, L. Paris, \textit{Gaussian groups and Garside groups, two generalizations of Artin groups}, Proc. London Math. Soc. 79 (3), (1999), 569-604.
\bibitem{em} E. ElRifai, H. Morton, \textit{Algorithms for positive braids}, Quart. J. Math. Oxford. Ser. (2) 45 (1994), no. 180, 479-497.
\bibitem{thurston} D. B. A. Epstein, J. Cannon, D. Holt, S. Levy, M. Paterson, W. Thurston, \textit{Word processing in groups}, Jones and Bartlett Publishers, Boston, MA, 1992.
\bibitem{Farb} B. Farb, D. Margalit, \textit{A primer on Mapping Class Groups}, Princeton Mathematical Series, 2011. 
\bibitem{Fathi} A. Fathi, F. Laudenbach, V. Poenaru, \textit{Travaux de Thurston sur les surfaces}, Ast\'{e}risque 66-67, SMF 1991/1979.
\bibitem{garside} F. Garside, \textit{The braid groups and other groups}, Quart. J. Math. Oxford Ser. (2) 20 (1969), 235-254.
\bibitem{Gebhardt} V. Gebhardt, \textit{A new approach to the conjugacy problem in Garside groups}, J. Algebra 292 (1), (2005), 282-302.
\bibitem{GebhardtGM} V. Gebhardt, J. Gonz\'{a}lez-Meneses, \textit{The cyclic sliding operation in Garside groups}, Math. Z. 265 (1), (2010), 85-114.
\bibitem{GebhardtGMalgorithm} V. Gebhardt, J. Gonz\'{a}lez-Meneses, \textit{Solving the conjugacy problem in Garside groups by cyclic sliding}, Journal of Symbolic Computation 45 (6) (2010), 629-656.
\bibitem{juanjaca} J. Gonz\'{a}lez-Meneses, \textit{On reduction curves and Garside properties of braids}, Contemp. Math. 538 (2011), 227-244.
\bibitem{GMWiest} J. Gonz\'{a}lez-Meneses, B. Wiest, \textit{Reducible braids and Garside theory}, Algebraic and Geometric Topology 11 (2011), 2971-3010. 
\bibitem{leelee} E.-K. Lee, S.-J. Lee, \textit{A Garside-theoretic approach to the reducibility problem in braid groups}, J. Algebra 320 (2008), no. 2, 783-820.
\bibitem{MasurMinsky} H. Masur, Y. Minsky, \textit{Geometry of the complex of curves. II. Hierarchical structure}, Geom. Funct. Anal. 10 (2000), no. 4, 902-974.
\bibitem{prasolov} M.V. Prasolov, \textit{Small braids with large ultra summit set}, Mat. Zametki 89 (4) (2011), 577-588.
\bibitem{jintao} J. Tao, \textit{Linearly Bounded Conjugator Property for Mapping Class Groups}, Geom. Funct. Anal. 23 (2013), no. 1, 415-466.
\end{thebibliography}
\end{document}